\documentclass[11pt]{amsart}
\usepackage{amssymb, amsfonts, amscd}
\input xy
\xyoption{all}
\renewcommand{\tilde}{\widetilde}
\newcommand{\liminv}[1]{{\displaystyle{\mathop{\rm 
lim}_{\buildrel\longleftarrow\over{#1}}}}\,}

\newcommand{\bbC}{{\mathbb C}}
\newcommand{\C}{{\mathbb C}}
\newcommand{\bbZ}{{\mathbb Z}}
\newcommand{\Z}{{\mathbb Z}}

\newcommand{\bbP}{{\mathbb P}}

\renewcommand{\P}{{\mathbb P}}
\newcommand{\Oh}{{\mathcal O}}
\newcommand{\frakX}{{\mathfrak X}}
\renewcommand{\bar}{\overline}
\renewcommand{\hat}{\widehat}
\DeclareFontFamily{OT1}{rsfs}{}
\DeclareFontShape{OT1}{rsfs}{n}{it}{<->rsfs10}{}
\DeclareMathAlphabet{\script}{OT1}{rsfs}{n}{it}
\DeclareMathOperator{\VV}{V}

\DeclareMathOperator{\HH}{H}

\DeclareMathOperator{\RR}{R}
\DeclareMathOperator{\EE}{E}
\DeclareMathOperator{\Cl}{Cl}
\DeclareMathOperator{\Alb}{Alb}
\DeclareMathOperator{\Ext}{Ext}
\DeclareMathOperator{\Pic}{Pic}
\DeclareMathOperator{\Div}{Div}
\DeclareMathOperator{\Lef}{Lef}

\DeclareMathOperator{\ALeff}{ALeff}
\DeclareMathOperator{\im}{Im}

\DeclareMathOperator{\coker}{Coker}

\newcommand{\scrL}{\script L}
\newcommand{\scrM}{\script M}

\newcommand{\scrN}{\script N}

\newcommand{\scrE}{\script E}
\newcommand{\scrF}{\script F}
\newcommand{\scrG}{\script G}

\newcommand{\onto}{\twoheadrightarrow}
\newcommand{\into}{\hookrightarrow}
\newcommand{\by}[1]{\xrightarrow{#1}}

\newcommand{\tensor}{\otimes}
\newcommand{\isom}{\cong}

\theoremstyle{plain}
\newtheorem{lemma}{Lemma}[section]

\newtheorem{thm}{Theorem}
\newtheorem{prop}[lemma]{Proposition}
\newtheorem{cor}[lemma]{Corollary}

\theoremstyle{definition}
\newtheorem{remark}[lemma]{Remark}

\newtheorem{defn}{Definition}

\newenvironment{diagram}[1]{\arraycolsep=\doublerulesep\begin{array}{#1}
}{\end{array}}

\begin{document}
\date{}
\title{The Grothendieck-Lefschetz theorem for Normal Projective Varieties}
\author{G.~V.~Ravindra}
\address{Department of Mathematics, Washington
University, St. Louis, MO 63130, USA}
\email{ravindra@math.wustl.edu}
\author{V.~Srinivas}
\address{School of Mathematics, Tata Institute of
Fundamental Research, Homi Bhabha Road, Mumbai-400005,
India}
\email{srinivas@math.tifr.res.in}
\begin{abstract}
We prove that for a normal projective variety $X$ in
characteristic 0, and a base-point free ample line bundle $L$ on it, the
restriction map of divisor class groups  $\Cl(X)\to \Cl(Y)$ is an
isomorphism for a general member $Y\in |L|$  provided that $\dim{X}\geq
4$. This is a generalization of the Grothendieck-Lefschetz Theorem, for
divisor class groups of singular varieties.
\end{abstract}
\thanks{Srinivas was partially supported by a
Swarnajayanthi Fellowship of the D.S.T.}
\maketitle
We work over $k$, an algebraically closed field of characteristic $0$.

Let $X$ be a smooth projective variety over $k$ and $Y$ a smooth
complete intersection subvariety of $X$. The Grothendieck-Lefschetz
theorem states that if dimension $Y \geq 3$, the Picard groups of
$X$ and $Y$ are isomorphic.

In this paper, we wish to prove an analogous statement for singular
varieties, with the Picard group replaced by the divisor class group. 

Let $X$ be an irreducible projective variety which is regular in
codimension 1 (for example, $X$ may be irreducible and normal). Recall
that for such $X$, the divisor class group $\Cl(X)$
is defined as the group of linear equivalence classes of Weil divisors on
$X$ (see \cite{Ha}, II, \S6). If $\dim X=d$, then
$\Cl(X)$ coincides with the Chow group $CH_{d-1}(X)$ as defined in
Fulton's book \cite{Fulton}. If $Y\subset X$ is an irreducible Cartier
divisor, which is also regular in codimension 1, there is a well-defined
restriction homomorphism\footnote{The terminology is from the
non-singular case, where one is considering restriction of line bundles.}
\[\Cl(X)\to \Cl(Y),\] 
determined by the rule
\[[D]\mapsto [D\cap Y],\]
where $D$ is any irreducible Weil divisor in $X$ distinct from $Y$, and
$[D\cap Y]$ denotes the Weil divisor on $Y$ associated to the
intersection scheme $D\cap Y$. This may be viewed as a
particular case of the refined Gysin homomorphism $CH_i(X)\to
CH_{i-1}(Y)$ defined in \cite{Fulton}, for $i=\dim X-1$. 

Now let $X$ be an irreducible projective variety over $k$, regular in
codimension 1, and let $\scrL$ be an ample line bundle over $X$, together
with a linear subspace $\VV\subset \HH^{0}(X, \scrL)$ which gives a
base point free ample linear system $|{\VV}|$ on $X$. Let $Y\in
|{\VV}|$ be a general element of this linear system; by Bertini's
theorem, we have $Y_{sing}=Y\cap X_{sing}$. In this context, our main
result is the following, which is an analogue of the 
Grothendieck-Lefschetz theorem.
\begin{thm}\label{thm1}
In the above situation, for a dense Zariski open set of $Y\in |{\VV}|$,
the restriction map
\[\Cl(X)\to \Cl(Y)\]
is an isomorphism, if $\dim X\geq 4$, and is
injective, with finitely generated cokernel, if $\dim
X=3$.
\end{thm}

Our proof is purely algebraic, in the style of the proof of the
Grothendieck-Lefschetz theorem given in \cite{H}, Chapter~IV. The above
result has an application in the theory of Deligne's 1-motives
(see \cite{D}), which is discussed in \S4 below; for this, it is of
interest to have such an algebraic proof. In an appendix, we also
sketch a different, transcendental proof of the theorem, when $k=\C$, due
to N. Fakhruddin, using results from stratified Morse theory, and
properties of the weight filtration on cohomology. 

We would like to thank the referee for a critical reading of the paper,
and detailed suggestions for improvement, including the discussion in \S5. 

\begin{section}{The Grothendieck-Lefschetz theorem for big linear systems}

In the situation of Theorem~\ref{thm1}, if $\tilde X \by{\pi} X$ is a
desingularisation of $X$, we have the following (Cartesian) diagram:
\[
\begin{diagram}{ccc}
\tilde Y & \into & \tilde X \\
\downarrow & & \downarrow \\
Y & \into & X \\
\end{diagram}
\]
Note that $\tilde Y$ is a general member of the pull-back linear
system $\pi^{\ast}\VV$ on the smooth proper variety $\tilde X$, and
therefore is smooth, by Bertini's theorem; hence $\tilde Y \to Y$ is a
desingularisation of $Y$.  If $X$ is singular, then $\tilde Y$ is a
general member of the linear system determined by 
$\pi^{\ast}V\subset \HH^{0}(\tilde X, \pi^{\ast}\scrL)$ where
$\pi^{\ast}\scrL$ is not ample, but is \emph{big and base-point free}.

Let $E_X=\pi^{-1}(X_{sing})\subset \tilde{X}$ be the exceptional
locus. Then $E_Y=E_X\cap\tilde{Y}$ is the exceptional locus for
$\tilde{Y}\to Y$. We have natural isomorphisms 
\[\Cl(X)\cong\Cl(X\setminus X_{sing})\cong\Pic(X\setminus X_{sing})\cong
\Pic(\tilde{X}\setminus E_X),\]
\[\Cl(Y)\cong\Cl(Y\setminus Y_{sing})\cong\Pic(Y\setminus Y_{sing})\cong
\Pic(\tilde{Y}\setminus E_Y).\]
This is because (i) the divisor class group is unchanged upon removal of a
closed subset of codimension $\geq 2$, and (ii) the divisor class group
coincides with the Picard group, for non-singular varieties. Thus,
Theorem~\ref{thm1} may be viewed as an assertion about the natural
restriction homomophism 
\[\Pic(\tilde{X}\setminus E_X)\to\Pic(\tilde{Y}\setminus E_Y).\] 
Also, by a standard argument (repeated below),  
$\Pic(\tilde{X})\to\Pic(\tilde{X}\setminus E_X)$ is surjective, with
kernel isomorphic to the free abelian group on the irreducible
divisorial components of $E_X$, and a similar assertion holds for
$\Pic(\tilde{Y})\to \Pic(\tilde{Y}\setminus E_Y)$. Indeed, since
$\tilde{X}$, $\tilde{Y}$ are non-singular, any line bundle on any Zariski
open subset extends to a line bundle on the variety, so the two
restriction maps are surjective, with kernel given by the line bundles
associated to divisors with support in $E_X$ and $E_Y$ 
respectively. However, if $E$ is any non-zero divisor
on $\tilde{X}$ with support in $E_X$, then $\Oh_{\tilde{X}}(E)$ is
a non-trivial line bundle: if not, we would have a non-constant rational
function $f$ on $\tilde{X}$ with divisor $E$; then $f$ determines a
non-zero regular function on $X\setminus X_{sing}$, which must extend to a
regular function on the normalization $X_n$ of $X$. But $X_n$ is an
irreducible projective variety, so any global regular function on it is
constant, which is a contradiction. This argument applies to
$\tilde{Y}$ as well. 

Thus, Theorem~\ref{thm1} is a consequence 
of the following version of the Grothendieck-Lefschetz
theorem for a big and base-point free linear system,
describing the kernel and cokernel of the restiction
map on Picard groups, for the inclusion
$\tilde{Y}\into \tilde{X}$ of a general member of the
linear system. The statement is a perhaps a bit technical, but
there is an obvious geometric motivation for the conditions stated.

\begin{thm}~\label{thm1'}
Let $\tilde{X}$ be a non-singular projective
$k$-variety, $\scrM$ a big invertible sheaf,
$\VV\subset H^0(\tilde{X},\scrM)$ a $k$-subspace
giving a base-point free linear system on $\tilde{X}$.
Let $\varphi:\tilde{X}\to\P^N_k$ be the morphism
determined by $|\VV|$, and $\tilde{X}\by{\pi}X\to
\P^N_k$ be the Stein factorization of $\varphi$.
Suppose $\dim \tilde{X}\geq 3$. Then for a Zariski
open subset of divisors $\tilde{Y}\in |\VV|$, the
restriction map
\[\rho:\Pic(\tilde{X})\to\Pic(\tilde{Y})\] has the
following properties.
\begin{enumerate}
\item[(a)]
$\rho$ has kernel (freely) generated by the classes of the irreducible
divisors  $E\subset\tilde{X}$ with $\dim\pi(E)=0$, and $\rho$ has a
finitely generated cokernel.
\item[(b)] If $F$ is a divisor on $\tilde{Y}$ supported in $E_Y$, such
that $\Oh_{\tilde{Y}}(F)\in {\rm image}\,\Pic(\tilde{X})$, then there is a
divisor $E$ on $\tilde{X}$ supported in $E_X$ with
$\Oh_{\tilde{X}}(E)\tensor\Oh_{\tilde{Y}}\cong
\Oh_{\tilde{Y}}(F)$. 
\item[(c)] If $\dim \tilde{X}\geq 4$, then the classes
of $\Oh_{\tilde{Y}}(E)$, with $E$ supported in $E_Y$, generate
$\coker(\rho)$. \end{enumerate} \end{thm}

\section{Some lemmas on vanishing of cohomology}

Next, we collect together some technical lemmas, used in the proof of 
Theorem~\ref{thm1'}.  

We first state a lemma due to Grothendieck (see
\cite{AG2,Mehta-Srinivas}). 
\begin{lemma}[Artin-Rees formula]\label{artinrees}
Let $f:V \to W$ be a proper morphism of Noetherian schemes, $\scrF$ a
coherent sheaf on $V$, $\mathcal I \subset \Oh_{W}$ a coherent ideal
sheaf, and $\mathcal J = f^{-1}\mathcal I.\Oh_{W}$ the inverse image
ideal sheaf on $V$. There exists $n_{0} \ge 0$ such that for all $n
\ge n_{0}$, the natural map 
$$\mathcal I^{n-n_{0}}\tensor\RR^{i}f_{\ast}(\mathcal J^{n_{0}}\scrF)
\to \RR^{i}f_{\ast}(\mathcal J^{n}\scrF)$$ is a surjection.
\end{lemma}
Here, $\mathcal J^n\mathcal F$ denotes the image of the
multiplication map $\mathcal J^n\tensor\mathcal F\to \mathcal F$.

\begin{lemma}\label{useful1}
Let $\tilde W \by{\pi} W$ be a proper surjective morphism, where
$\tilde{W}$ is 
an irreducible non-singular variety of dimension
$d$, and $\dim W\ge 2$. Let $\mathcal F$ be a coherent sheaf on $\tilde
W$, with $R^{d-1}\pi_*\mathcal F\neq 0$. Then there exists an effective
divisor $E\subset \tilde W$ whose support has 0-dimensional image under
$\pi$, such that $\RR^{d-1}\pi_{\ast}\mathcal F(-E) \to 
\RR^{d-1}\pi_{\ast}\mathcal F$ is
the zero map.
\end{lemma}
\begin{proof}
Let $S\subset W$ be the support of $\RR^{d-1}\pi_*\mathcal F$. Then $S$
consists of points $w\in W$ with $\dim\pi^{-1}(w)\geq d-1$, and under our
hypotheses, this forces $\dim S=0$. 

Let $\mathcal I\subset\Oh_W$ be the ideal sheaf of $S$, 
and let 
\[\mathcal J={\rm image}\,
\left(\pi^*\mathcal I\to \Oh_{\tilde W}\right)\]
be the inverse image ideal sheaf in $\Oh_{\tilde W}$. 
Lemma~\ref{artinrees} above implies that there exists an  $m_{0} \ge 0$
such that the map
$$\mathcal I^{m-m_{0}}\tensor\RR^{d-1}\pi_{\ast}(\mathcal J^{m_0}\mathcal
F) \to \RR^{d-1}\pi_{\ast}(\mathcal J^m \mathcal F)$$ 
is a surjection.

We claim that, for large enough $m$, the map
$\RR^{d-1}\pi_{\ast}(\mathcal J^m\mathcal F) \to
\RR^{d-1}\pi_{\ast}(\mathcal F)$ is the zero map. This is because by
Lemma~\ref{artinrees}, we have a diagram
\[
\xymatrix{ \mathcal
I^{m-m_{0}}\tensor\RR^{d-1}\pi_{\ast}
(\mathcal J^{m_0}\mathcal F) \ar@{-{>>}}[r]
\ar[d] &
\RR^{d-1}\pi_*(\mathcal J^m\mathcal F)\ar[d]\\ 
\mathcal I^{m-m_0}\tensor \RR^{d-1}\pi_{\ast}(\mathcal F) \ar[r] &
\RR^{d-1}\pi_{\ast}(\mathcal F) \\}
\]
where the top horizontal arrow is surjective. The bottom horizontal map is
0, if $m$ is large enough, since $\mathcal I$ is the ideal
defining the support of $\RR^{d-1}\pi_*(\mathcal F)$; hence the right
vertical arrow is 0. 

Since $\tilde W$ is non-singular, there exists an effective
(Cartier) divisor $E_0$ in $\tilde
W$ such that $\mathcal J=\Oh_{\tilde W}(-E_0) \tensor J$, where 
$J\subset \Oh_{\tilde W}$ defines a subscheme of codimension $\ge 2$. 
In particular, we have inclusions of ideal sheaves $\mathcal J^m\subset
\Oh_{\tilde W}(-mE_0)\subset \Oh_{\tilde W}$.

 Now consider the exact
sequence
\[
0 \to \mathcal J^m\mathcal F \to \mathcal
F(-mE_0) \to \mathcal F(-mE_0)\tensor\Oh_{mZ} \to 0
\]
Here $mZ\subset \tilde W$ is the subscheme defined by the ideal sheaf $J^m
\subset \Oh_{\tilde W}$.  This gives a long exact sequence
\[
\RR^{d-1}\pi_{\ast}(\mathcal J^m\mathcal F) \to
 \RR^{d-1}\pi_{\ast}(\mathcal F(-mE_0))\to \RR^{d-1}\pi_{\ast}(\mathcal
 F(-mE)_{\vert{mZ}})
\]
We note that the last term is zero since the codimension of $Z$ in
$\tilde W$ is $\ge 2$. Hence the first arrow is a surjection.

To conclude the proof of the lemma, we note that the (zero) map 
$\RR^{d-1}\pi_{\ast}(\mathcal J^m\mathcal F) \to 
\RR^{d-1}\pi_{\ast}(\mathcal F)$
factors as 
$$\RR^{d-1}\pi_{\ast}({\mathcal J^m}\mathcal F) \to
\RR^{d-1}\pi_{\ast}(\mathcal F(-mE_0))\to \RR^{d-1}\pi_{\ast}(\mathcal
F).$$ 
Since the first map is a surjection, the second is necessarily
the zero map. The lemma thus holds with $E=mE_0$, where $m$ is
sufficiently large.
\end{proof}

Another version of the above lemma is available, in a situation analogous
to Lemma~\ref{redn}.  
\begin{lemma}\label{useful2}
Let $\tilde W \by{\pi} W$ be a desingularisation of a normal
projective variety of dimension $d$ and $\mathcal F$ be a coherent
sheaf on $\tilde W$. Assume that there exists an effective divisor $E$ 
on $\tilde{W}$ with $\pi$-exceptional support, such
that $-E$ is $\pi$-ample. 
\begin{enumerate}
\item[(a)] There exists a positive integer $r_0$ such that 
such that $\RR^i\pi_{\ast}\mathcal F(-rE)=0$ for all $r\geq r_0$ and
all $i>0$.  
\item[(b)] Suppose $\mathcal L$ is ample on $W$, and $\mathcal F$ is 
locally free on $\tilde{W}$. Then there exists a positive integer $r_1$
such that for each $r\geq r_1$, and for all $n\geq n_1$ (depending on
$r$), we have
\[H^i(\tilde{W},\mathcal F(rE)\tensor\pi^*\mathcal
L^{-n})=0\;\;\;\mbox{ for all $i< d$},\]
\[H^i(\tilde{W},\mathcal F(-rE)\tensor\pi^*\mathcal
L^{n})=0\;\;\;\mbox{ for all $i>0$}.\] 
\end{enumerate}
\end{lemma}
\begin{proof} The assertion in (a) is just Serre's vanishing theorem,
since $\Oh_{\tilde{W}}(-E)$ is $\pi$-ample. The two assertions in (b) 
are equivalent (with perhaps different values of $r_1$, $n_1$), using
Serre duality on $\tilde{W}$. The second assertion in (b) follows from
(a), using the Leray spectral sequence for $\pi$, and Serre's vanishing on
$W$ for the ample line bundle $\mathcal L$.
\end{proof}

We recall a form of the Grauert-Riemenschneider theorem which we need
below.
\begin{thm}\label{GR}
Let $W$ be a non-singular projective variety over $k$, and $\scrM$ a big
and base-point free line bundle on $W$. Then $H^i(W,\scrM^{-1})=0$ for
$i<\dim W$.
\end{thm}
\begin{proof} 
A proof of this statement can be found in \cite{EV}, Cor. 5.6(b).
\end{proof}

\section{Proof of Theorem~2}
\subsection{Some preliminary reductions} 
From now on, we fix the following notation (used in
Theorem~\ref{thm1'})  -- 
$\tilde{X}$ is a smooth projective $k$-variety, $\mathcal M$ a big line
bundle on $\tilde{X}$, $\VV\subset H^0(\tilde{X},\mathcal M)$ a subspace
giving a linear system $|\VV|$  without base points, $\tilde{Y}$ a general
member of this linear system, $\pi:\tilde X\to X$ obtained by Stein
factorization of the morphism determined by $|\VV|$, $\tilde Y\to Y$ the
induced morphism (which is also the Stein factorization of the restriction
to $\tilde Y$ of the original morphism on $\tilde X$), and $\scrL$ the
invertible sheaf on $X$ such that $\pi^*\scrL=\scrM$. 

We first make a further reduction. 
\begin{lemma}\label{redn} To prove
Theorem~\ref{thm1'}, it suffices to do it in the case
when the morphism $\pi:\tilde{X}\to X$ (obtained by
Stein factorization) has a purely divisorial
exceptional locus, with non-singular irreducible components, and there
exists a $\pi$-ample divisor of the form $-E$ where $E$ is an effective
divisor with $\pi$-exceptional support.
\end{lemma}
\begin{proof}  Since $\tilde{X}\to X$ and $\tilde{Y}\to
Y$ are obtained from Stein factorizations of the morphisms determined
by the base point free linear system $|\VV|$, we have that  $X$, $Y$ are
normal projective varieties, such that $Y$ is a Cartier divisor in $X$.
There is an induced restriction homomorphism $\Cl(X)\to\Cl(Y)$. 

It is easy to see that the conclusions of Theorem~\ref{thm1'} hold for
$\Pic\tilde{X}\to\Pic\tilde{Y}$ if and only if the
conclusions of Theorem~\ref{thm1} hold for 
$\Cl(X)\to\Cl(Y)$ (i.e., the restriction map
on class groups is an isomorphism if $\dim X\geq
4$, and an inclusion with finitely generated
cokernel if $\dim X=3$). 

In particular, if we replace $\tilde{X}\to X$ by
another resolution of singularities $\pi':\tilde{X}'\to
X$, and $\tilde{Y}$ by the inverse image
$\tilde{Y}'$ of $Y$ in that resolution, it suffices
to prove Theorem~\ref{thm1'} for this new pair
$(\tilde{X}',\tilde{Y}')$, and the pull-back linear 
system from $X$ (note that there is an open subset of 
the linear system consisting of divisors which are
``general'' for both $\tilde{X}$ and $\tilde{X}'$).  

Now by Hironaka's theorem, we can find a resolution of singularities
$\pi':\tilde{X}'\to X$ such that the exceptional locus is divisorial,
with non-singular irreducible components. Further,
there is an effective exceptional divisor $E$ such
that $-E$ is $\pi'$-ample; this is because the resolution may be
obtained by successive blow-ups at centres lying over the singular
locus, so that the resolution may be viewed as a blow-up of an ideal
sheaf whose radical defines the singular locus. The pull-back ideal
sheaf is an invertible sheaf which is $\pi'$-ample,
and is the ideal sheaf on $\tilde{X}'$ of some subscheme with exceptional
support.
\end{proof}
\begin{remark}
This reduction is needed only in the proof that, if $\dim X\geq 4$, then
$\Cl(X)\to\Cl(Y)$ is surjective. 
\end{remark}

We follow the line of proof of the Grothendieck-Lefschetz theorem as
given in \cite{H}. The idea is to pass from $\tilde Y$ to the formal
completion $\frakX$ of $\tilde X$ along $\tilde Y$. Then from
$\frakX$ we pass to a neighbourhood $U$ of $\tilde Y$ in
$\tilde X$, using a version of the Lefschetz Conditions, and then to
$\tilde X$ itself.

\begin{lemma}\label{exp}
With notation as above, if $\dim \tilde{X}\geq 4$, then $\Pic(\frakX) 
\cong \Pic(\tilde Y)$. If $\dim\tilde{X}=3$, then
$\Pic(\frakX)\to\Pic(\tilde{Y})$ is injective, with finitely
generated cokernel.
\end{lemma}

\begin{proof}
Consider the short exact sequence
\[
0 \to \scrM^{-m}\tensor\Oh_{\tilde Y} \to
\Oh^{\times}_{\frakX_{m+1}} \to \Oh^{\times}_{\frakX_{m}}
\to 0
\]
where $\frakX_{m}$ is the $m$-th infinitesimal neighbourhood of
$\tilde Y \subset \tilde X$ and in particular $\frakX_{1}=\tilde
Y$. As usual, $\Oh^{\times}_T$ denotes the (multiplicative) sheaf of
invertible regular  functions on $T$. The first horizontal sheaf map is
the ``exponential map'', defined on sections by $s\mapsto 1+s$.  Taking
the cohomology long exact sequence, one has 
\[
\begin{diagram}{ccccccccc}
 {\scriptstyle \to} & {\scriptstyle\HH^{1}(\tilde Y,
\scrM^{-m}\tensor\Oh_{\tilde Y})} & {\scriptstyle\to} &
{\scriptstyle \HH^{1}(\tilde \frakX_{m+1},\Oh_{
\frakX_{m+1}}^{\times})} & {\scriptstyle\to} &
{\scriptstyle\HH^{1}( \frakX_{m}, \Oh_{\frakX_{m}}^{\times})}
& {\scriptstyle \to} & {\scriptstyle\HH^{2}(\tilde
Y,\scrM^{-m}\tensor\Oh_{\tilde Y})} & {\scriptstyle\to}
\end{diagram}
\]

By the Grauert-Riemenschneider vanishing theorem (Theorem~\ref{GR} above),
the extreme terms vanish if $\dim \tilde{Y}\geq 3$, and thus we have
$\Pic(\frakX_{m+1}) \isom \Pic(\frakX_{m})$ for each $m\geq 1$.  
From the Grothendieck formula (see \cite{Ha} II Ex. 9.6, for example)
\[\Pic(\frakX)\cong \liminv{m}\Pic(\frakX_m),\]
we then get $\Pic(\frakX) \isom \Pic(\tilde Y)$.

If $\dim\tilde{Y}=2$, then the same argument shows that
$\Pic(\frakX_m)\to\Pic(\tilde{Y})$ is injective for each $m$. On
the other hand, $H^1(\tilde{X},\Oh_{\tilde{X}})\to
H^1(\tilde{Y},\Oh_{\tilde{Y}})$ is an isomorphism, since
$H^i(\tilde{X},\Oh_{\tilde{X}}(-\tilde{Y}))$ vanishes for $i<3$
(Grauert-Riemenschneider vanishing, Theorem~\ref{GR} above). Hence
$\Pic^0(\tilde{X})\to \Pic^0(\tilde{Y})$ is an isogeny, and in particular
is surjective. Hence $\coker\Pic(\tilde{X})\to\Pic(\tilde{Y})$ is a
quotient of the Neron-Severi group of $\tilde{Y}$, and is finitely
generated. A similar conclusion then clearly holds for
$\coker\Pic(\frakX)\to\Pic(\tilde{Y})$. 
\end{proof}

\begin{subsection}{The condition $\Lef(\tilde X, \tilde Y)$}

In the proof of the Grothendieck-Lefschetz theorem (see \cite{H}, Ch. IV),
one considers the \emph{Lefschetz condition} which implies the injectivity
of the morphism between the Picard groups. We will show that it holds in
our situation as well.

If $\scrE$ is a coherent sheaf on some open neighbourhood of $\tilde{Y}$
in $\tilde{X}$, then $\widehat\scrE$ denotes the corresponding
(formal) coherent sheaf on the formal completion $\frakX$ of
$\tilde{X}$ along $\tilde{Y}$. With this notation, recall the following
definition (see \cite{H}, page~164; this is a slight modification of
Grothendieck's definition in \cite{AG2}, page~112, as remarked by the
referee).
\begin{defn}
The pair $(\tilde X, \tilde Y)$ {\em satisfies the Lefschetz condition
$\Lef(\tilde X, \tilde Y)$} if for every open set $U\subset \tilde{X}$
containing $\tilde Y$, and every locally free sheaf $\scrE$ on
$U$, there exists an open set $U'$ with $\tilde Y
\subset U' \subset U$ such that the natural map
$$ \HH^{0}(U',\scrE_{|U'}) \to \HH^{0}(\frakX,
\widehat\scrE)$$ is an isomorphism.
\end{defn}

Note that there is a finite (perhaps empty) set $S_0\subset X$ of
(closed) points $x\in X$ with $\dim \pi^{-1}(x)=\dim X-1$. Since $\tilde
Y$ is general, we may assume that $Y\cap S_{0}=\emptyset$. Further note
that any divisor in $\tilde X$ whose support is disjoint from $\tilde Y$
must be supported in  $\pi^{-1}(S_{0})$. If $E$ is any divisor on
$\tilde{X}$ with support in $\pi^{-1}(S_{0})$, then
$\widehat{\Oh_{\tilde{X}}(E)}\cong\Oh_{\frakX}$.

Recall that the dual of a coherent sheaf $\scrN$ on $\tilde{X}$ is
$\scrN^{\vee}={\mathcal Hom}_{\Oh_{\tilde{X}}}(\scrN,\Oh_{\tilde{X}})$; 
recall also that $\scrN$ is {\em reflexive} if
$\scrN\to(\scrN^{\vee})^{\vee}$ is an isomorphism.

\begin{lemma}\label{alef}
Let $\scrN$ be a reflexive, coherent sheaf on $\tilde
X$ which is locally free in a neighbourhood of $\tilde{Y}$. Then there
exists an effective divisor $E$ on $\tilde{X}$, where either $E=0$ or 
$\dim \pi({\rm supp}\,E) =0$,  such that the natural map
$$\HH^{0}(\tilde X, \scrN(E)) \to
\HH^{0}(\frakX,\widehat\scrN) $$
is an isomorphism. 
\end{lemma}
\begin{proof}
Let $n=\dim\tilde{X}$.  Using Serre duality on $\tilde X$ and formal
duality on $\frakX$ (see \cite{H}, III, Theorem~3.3), 
 we reduce to proving that
$$ \HH^{n}_{\tilde Y}(\tilde X, \scrN^{\vee}(-E)\tensor \omega_{\tilde X})
\to \HH^{n}(\tilde X, \scrN^{\vee}(-E)\tensor \omega_{\tilde X})$$ 
is an isomorphism, for appropriate $E$. Here, we note that though $\scrN$
may not be locally free, Serre duality implies that the dual of the
finite dimensional vector space $H^n(\tilde{X},\scrN^{\vee}(-E)\tensor
\omega_{\tilde{X}})$ is 
\[{\rm Hom}(\scrN^{\vee}(-E)\tensor 
\omega_{\tilde{X}},\omega_{\tilde{X}})=H^0(\tilde{X},\scrN(E)),\]
since $\scrN$ is reflexive. 

For any effective divisor $E$ supported in $\pi^{-1}(S_0)$, consider the
commutative diagram with exact rows
\[
\begin{diagram}{cccccccc}
\HH^{n-1}(\tilde X \setminus \tilde Y,
\scrN^{\vee}(-E)\tensor
\omega_{\tilde X}) &  \by{\delta_{1}} & 
\HH^{n}_{\tilde Y}(\tilde X,
\scrN^{\vee}(-E)\tensor \omega_{\tilde X}) &  \onto
&\HH^{n}(\tilde X,
\scrN^{\vee}(-E)\tensor \omega_{\tilde X})\\
 \downarrow{\phi_{1}} & & \downarrow{\phi_{2}}
 & & \downarrow{\phi_{3}}
& & \\ 
\HH^{n-1}(\tilde X\setminus \tilde Y, \scrN^{\vee}\tensor
\omega_{\tilde X}) & {\scriptstyle \by{\delta_{2}}} & 
\HH^{n}_{\tilde Y}(\tilde X,
\scrN^{\vee}\tensor \omega_{\tilde X})) & \onto &
\HH^{n}(\tilde X,
\scrN^{\vee}\tensor \omega_{\tilde X})\\
\end{diagram}
\]

The last maps in the sequences are onto since $\tilde X \setminus
\tilde Y$ has cohomological dimension at most $n-1$. Moreover, the map
$\phi_{2}$ is an isomorphism by excision since $\tilde Y\cap E =
\emptyset$. The Leray spectral sequence for the map $\pi: \tilde X
\setminus \tilde Y \to X \setminus Y$ applied to the cohomology group
$\HH^{n-1}(\tilde X \setminus \tilde Y, \scrN^{\vee}(-E)\tensor
\omega_{\tilde X})$ has $\EE_{2}^{p,q}=0$ for $p>0$ ($X \setminus Y$ is
affine!). By Lemma~\ref{useful1}, there exists an $E$ as in the statement
of the lemma such that the map 
$\RR^{n-1}\pi_{\ast}(\scrN^{\vee}(-E)\tensor \omega_{\tilde X}) \to
\RR^{n-1}\pi_{\ast}(\scrN^{\vee}\tensor \omega_{\tilde X})$ is the zero
map, and thus the map $\phi_1$ is the zero map, for this choice of
$E$. This in turn implies that the corresponding map $\delta_1$ is zero.

We thus have the following commutative diagram
\[
\begin{diagram}{ccc}
\HH^{n}_{\tilde Y}(\tilde X,\scrN^{\vee}(-mE)\tensor \omega_{\tilde X}) &
\isom & \HH^{n}(\tilde X ,\scrN^{\vee}(-mE)\tensor \omega_{\tilde X}) \\
\downarrow{\isom} & & \downarrow \\ 
\HH^{n}_{\tilde Y}(\tilde X,\scrN^{\vee}\tensor \omega_{\tilde X}) & \onto
&
\HH^{n}(\tilde X ,\scrN^{\vee}\tensor \omega_{\tilde X}) \\
\end{diagram}
\]
Dualising, we have 
$$\HH^{0}(\tilde  X, \scrN(mE)) \isom
\HH^{0}(\frakX,\widehat\scrN)$$
\end{proof}
\begin{cor}
The condition $\Lef(\tilde X, \tilde Y)$ holds.
\end{cor}
\begin{proof}
For any open set $U\supset \tilde Y$ in $\tilde X$, and any locally
free sheaf $\scrN_U$ on $U$, we can find a reflexive sheaf $\scrN$ on
$\tilde{X}$ extending $\scrN_U$, i.e., with $\scrN\mid_{U}\cong
\scrN_U$ (first choose a coherent extension, then replace it by its
double dual).  For a suitable divisor $E$ with $\dim({\rm
supp}\,E)=0$, we have a commutative diagram induced by restriction
maps
\[
\begin{diagram}{ccc}
\HH^{0}(\tilde X, \scrN(E)) & \isom &
\HH^{0}(\frakX,\widehat\scrN) \\ \searrow & & \nearrow \\ &
\HH^{0}(U, \scrN(E)) & \\
\end{diagram}
\]
In particular, for any open $V$ such that $\tilde Y \subset V
\subset U$ and $V\cap {\rm supp}\,E = \emptyset$ the above
factorisation gives a surjection
$$\HH^{0}(V, \scrN(E)) \isom \HH^{0}(V, \scrN)\onto
\HH^{0}(\frakX,\widehat\scrN).$$ 
But since $\scrN$ is locally free on $V$ and $V$ is irreducible, 
the map is also an injection. Thus $\Lef(\tilde X, \tilde Y)$ holds.
\end{proof}
\begin{cor}\label{Lef(X,Y)}
For normal $X$, $Y$ as above, the condition $\Lef(X,Y)$ holds.
\end{cor}
\begin{proof}
Since $X$ and $Y$ are normal, one has $\pi_{\ast}\Oh_{\tilde
X}\isom\Oh_{X}$ and $\pi_{\ast}\Oh_{\tilde Y} \isom
\Oh_{Y}$. $\Lef(X,Y)$ then follows from $\Lef(\tilde X,\tilde Y)$
applied to sheaves which are pull backs of locally free sheaves on
neighbourhoods of $Y$ in $X$. 
\end{proof}
\begin{cor}\label{kernel}
The kernel of the restriction map $\Pic(\tilde X) \to
\Pic(\tilde Y)$ is freely generated by the
classes of irreducible effective divisors which map to points in $X$.
\end{cor}
\begin{proof}

It is obvious that the classes of such divisors are contained in the
kernel since $Y\cap S_0=\emptyset$ (as $Y$ is general). On the other hand,
if $\mathcal N$ is a line bundle on $\tilde{X}$ with 
$\mathcal N\tensor\Oh_{\tilde Y}\cong\Oh_{\tilde Y}$, 
then we first note that
$\widehat{\mathcal N}\cong \Oh_{\frakX}$ by Lemma~\ref{exp}, and there is
thus an invertible element of $H^0(\frakX,\widehat{\mathcal N})$; by
Lemma~\ref{alef}, this formal global section is obtained from a global
section on $\tilde{X}$ of $\mathcal N(E)$ for some divisor $E$ on
$\tilde{X}$ supported over $S_0\subset X$. This section of $\mathcal N(E)$
has no zeroes when restricted to $\tilde{Y}$, so its divisor of zeroes
$E'$ is also supported over $S_0$, and hence $\mathcal N\cong\Oh_{\tilde
X}(E'-E)$.  
\end{proof}
\begin{cor}\label{formalsects} For each $n\geq 0$, and any effective
divisor $F$ on $\tilde{X}$ with $\pi$-exceptional support, the natural
maps
\[H^0(\tilde{X},{\mathcal M}^{\tensor n})\to 
H^0(\tilde{X},{\mathcal M}^{\tensor n}(F))\to
H^0(\frakX,\widehat{{\mathcal M}^{\tensor n}(F)})\] 
are isomorphisms. 
\end{cor}
\begin{proof} Since $\scrM=\pi^*\scrL$, where $\scrL$ is invertible on the
normal variety $X$, and $\pi_*\Oh_{\tilde X}=\Oh_X$, it follows that
$\scrM^{\tensor n}\into\scrM^{\tensor n}(E+F)$ is an isomorphism on global
sections for any effective divisor $E$ supported in $\pi^{-1}(S_0)$, and
any $n\geq 0$.
\end{proof} 
\end{subsection}

\begin{subsection}{The condition {\bf $\ALeff(\tilde X,\tilde Y)$}}

We now introduce a second condition $\ALeff(\tilde X,\tilde Y)$
(Almost Effective Lefschetz), which is a variation of Grothendieck's
Effective Lefschetz Condition (denoted by ``Leff'' in \cite{H}). In the
proof in \cite{H} of the Grothendieck-Lefschetz theorem, the condition
Leff is used to show the surjectivity of the restriction map between the
Picard groups; ALeff has a similar role here.

\begin{defn} 
We say the pair $(\tilde X, \tilde Y)$ satisfies the {\em ALeff condition}
if 
\begin{enumerate}
\item $\Lef(\tilde X, \tilde Y)$ holds and 
\item for any (formal) 
invertible  sheaf $\scrE$ on $\frakX$ there
exists an open set $U$ containing $\tilde{Y}$ and 
an invertible sheaf $\mathcal E$ on $U$, together with a map 
$\widehat{\mathcal {E}}\to \scrE$, which is an isomorphism outside the 
exceptional locus of $\pi:\tilde{Y}\to Y$.
\end{enumerate}
\end{defn}
Note that the formal scheme $\frakX$ is a ringed space with
underlying
topological space $\tilde{Y}$, so the second condition above is
meaningful.

For any formal coherent sheaf $\scrF$ on $\frakX$, we will make the
following abuses of notation: for any divisor $D$ on $\tilde{X}$,  
let $\scrF(D)$ denote the formal coherent sheaf  
$\scrF\tensor_{\Oh_{\frakX}}\widehat{\Oh_{\tilde{X}}(D)}$, 
and for any coherent $\mathcal G$ on $\tilde{X}$, let 
$\scrF\tensor\mathcal G$ denote the formal coherent sheaf 
$\scrF\tensor_{\Oh_{\frakX}}\widehat{\mathcal G}$.

\begin{prop}\label{prop1}
Let $\tilde X$, $\tilde Y$ be as in Lemma~\ref{redn}, with $\dim \tilde
X\geq 3$, and let $E$ be an effective divisor on $\tilde{X}$ with
exceptional support such that $-E$ is $\pi$-ample. Then for any formal
locally free sheaf $\scrF$ on $\frakX$, there exists $r>0$ such that
for any $m\ge 0$, if we set 
\[\scrG_m=\scrF(rE)\tensor\scrM^{\tensor m},\]
then for all $m>>0$,
\[\coker \left(H^0(\frakX,\scrG_m)\tensor_k\Oh_{\frakX}\to
\scrG_m\right)\]
is supported on $\tilde{Y}\cap E$.
\end{prop}
\begin{proof}
The proof is in several steps. Let
$\scrF_n=\scrF\tensor\Oh_{\frakX_{n}}$,
for $n\geq 1$, be the sequence of locally free sheaves (on the sequence of
schemes  $\frakX_{n}$) associated to the formal locally free sheaf
$\scrF$. We have exact sequences
\[0\to \scrF\tensor\scrM^{\tensor m-n}\to \scrF\tensor\scrM^{\tensor m}\to
\scrF_n\tensor\scrM^{\tensor m}\to 0\]
for each $m\in\Z$ and $n>0$, where the ideal sheaf of $\tilde{Y}$ in
$\Oh_{\tilde{X}}$ is identified with $\scrM^{-1}$.

\begin{lemma}\label{l1} Let $d=\dim \tilde{Y}$. There exists $r_0>0$ such
that, for each $r\geq r_0$, all $m>>0$ (depending on $r$), and all $i>0$,
we have 
\[H^i(\tilde{Y},\scrF_1\tensor\Oh_{\tilde{X}}(-rE)\tensor\scrM^{\tensor
m})=0,\]
\[H^{d-i}(\tilde{Y},\scrF_1\tensor\Oh_{\tilde{X}}(rE)\tensor\scrM^{\tensor
-m})=0,\]
\end{lemma}
\begin{proof} We have that $\scrM=\pi^*\scrL$ where $\scrL$ is ample on
$X$. Now apply Lemma~\ref{useful2}(b).
\end{proof}

\begin{lemma}\label{l2}
There exists $r_0>0$ so that, for any $r\geq r_0$ and any $m\in\mathbb Z$,
the vector space
\[H^1(\frakX,\scrF(rE)\tensor\scrM^{\tensor m})\]
is finite dimensional.
\end{lemma}
\begin{proof}
We have exact sheaf sequences
\[0\to \scrF_1(rE)\tensor\scrM^{\tensor m-n}\to \scrF_{n+1}(rE)\tensor
\scrM^{\tensor m}\to \scrF_n(rE)\tensor\scrM^{\tensor m}\to 0.\]
Since $\dim\tilde{Y}\geq 2$, we have for each given $m$ that
\[H^1(\tilde{Y},\scrF_1(rE)\tensor\scrM^{\tensor m-n})=0\]
provided $n>>m$, from Lemma~\ref{l1}; thus for $n>>m$,
\[H^1(\frakX_{n+1},\scrF_{n+1}(rE)\tensor\scrM^{\tensor m})\to
H^1(\frakX_{n},\scrF_{n}(rE)\tensor\scrM^{\tensor m})\] 
is injective. Hence, in the Grothendieck formula
\[H^1(\frakX,\scrF(rE)\tensor\scrM^{\tensor
m})=\liminv{n}H^1(\frakX_{n},\scrF_n(rE)\tensor\scrM^{\tensor m}),\]
the maps in the inverse system on the right are, for $n>>0$, injective
maps of finite dimensional $k$-vector spaces (the above inverse limit
formula holds, because the corresponding inverse system for $H^0$ is an
inverse system of finite dimensional $k$-vector spaces, hence
satisfies the Mittag-Leffler condition (ML)). 
Thus the inverse limit $H^1(\frakX,\scrF(rE)\tensor\scrM^{\tensor
m})$ is a finite dimensional vector space. 
\end{proof}

Now consider the exact sequences
\[0\to \scrF(-rE)\tensor\scrM^{\tensor m}\to
\scrF(-rE)\tensor\scrM^{\tensor
m+1}\to \scrF_1(-rE)\tensor \scrM^{\tensor m+1}\to 0.\]
\begin{lemma}\label{l3}
There exists $r_1>0$ so that, for each $r\geq r_1$, and all $m>>0$
(depending on $r$), the map
\[H^1(\frakX,\scrF(-rE)\tensor\scrM^{\tensor m})\to
H^1(\frakX,\scrF(-rE)\tensor\scrM^{\tensor m+1})\]
is surjective.
\end{lemma}
\begin{proof}
It suffices to see that
\[H^1(\tilde{Y},\scrF_1(-rE)\tensor\scrM^{\tensor m+1})=0.\]
This follows from Lemma~\ref{l1}.
\end{proof}

Fix an $r$ large enough so that the conclusions of Lemma~\ref{l2} and
\ref{l3} hold. Define 
\[V_m={\rm image}\,H^1(\frakX,\scrF(-rE)\tensor\scrM^{\tensor m}) 
\by{\beta_m} H^1(\frakX,\scrF(rE)\tensor\scrM^{\tensor m}),\]
where the map is induced by the natural inclusion
$\Oh_{\tilde{X}}(-rE)\to\Oh_{\tilde{X}}(rE)$, determined by the
tautological section of $\Oh_{\tilde{X}}(2rE)$.
From Lemma~\ref{l2}, $V_m$ is a finite dimensional vector space, and from
Lemma~\ref{l3}, the natural maps 
\[H^1(\frakX,\scrF(-rE)\tensor\scrM^{\tensor m})\to
H^1(\frakX,\scrF(-rE)\tensor\scrM^{\tensor m+1})\]
induce surjections $V_m\to V_{m+1}$, for all large enough $m$. 
Hence $V_m\to V_{m+1}$ is in fact an isomorphism, for all large enough
$m$. Consider the commutative diagram of formal sheaves with exact rows
\[\begin{array}{ccc}
0\to \scrF(-rE)\tensor\scrM^{\tensor m}\to &
\scrF(-rE)\tensor\scrM^{\tensor
m+1} &\to \scrF_1(-rE)\tensor \scrM^{\tensor m+1}\to 0\\
\downarrow&\downarrow&\downarrow\\
0\to \scrF(rE)\tensor\scrM^{\tensor m}\to &
\scrF(rE)\tensor\scrM^{\tensor
m+1} &\to \scrF_1(rE)\tensor \scrM^{\tensor m+1}\to 0.
\end{array}\]
The vertical arrows are induced by the natural inclusion
$\Oh_{\tilde{X}}(-rE)\to\Oh_{\tilde{X}}(rE)$ (as above, in defining
$V_m$). There is an induced cohomology diagram with exact rows

\small{\small{
$\hspace{-2cm}\begin{array}{rccl}
\HH^0(\frakX,\scrF(-rE)\tensor\scrM^{\tensor m+1})\to&
\HH^0(\tilde{Y},\scrF_1(-rE)\tensor\scrM^{\tensor m+1})&\to 
\HH^1(\frakX,\scrF(-rE)\tensor\scrM^{\tensor m})\to & 
\HH^1(\frakX,\scrF(-rE)\tensor\scrM^{\tensor m+1}) \\
\downarrow\quad\quad\quad\quad\quad
&\quad\downarrow\alpha&\quad\downarrow\beta_m&
\quad\quad\quad\quad\downarrow\beta_{m+1}\\
\HH^0(\frakX,\scrF(rE)\tensor\scrM^{\tensor m+1})\to&
\HH^0(\tilde{Y},\scrF_1(rE)\tensor\scrM^{\tensor m+1})&\to  
\HH^1(\frakX,\scrF(rE)\tensor\scrM^{\tensor m})\to & 
\HH^1(\frakX,\scrF(rE)\tensor\scrM^{\tensor m+1}) 
\end{array}
$}}

Since $V_m\to V_{m+1}$ is an isomorphism, we see that 
\[{\rm image}\,\left( H^0(\frakX,\scrF(rE)\tensor\scrM^{\tensor m+1})\to
H^0(\tilde{Y},\scrF_1(rE)\tensor\scrM^{\tensor m+1})\right)\]
contains the subspace
\[{\rm image}\,H^0(\tilde{Y},\scrF_1(-rE)\tensor\scrM^{\tensor m+1}).\]
Since $m>>0$, the global sections of the sheaf
$\scrF_1(-rE)\tensor\scrM^{\tensor m+1}$ generate it on
$\tilde{Y}\setminus E$ (since the direct image of this sheaf on $Y$ is
globally generated, and $\tilde{Y}\to Y$ is an isomorphism outside
$\tilde{Y}\cap E$). Hence, the natural map between coherent formal sheaves
\[H^0(\frakX,\scrF(rE)\tensor\scrM^{m+1})\tensor_k\Oh_{\frakX}\to
\scrF(rE)\tensor\scrM^{\tensor m+1}\]
restricts to a surjection on $\tilde{Y}\setminus E$. This proves
Proposition~\ref{prop1}
\end{proof}
\begin{remark} The referee has pointed out that, in
Proposition~\ref{prop1}, we can in fact get the stated conclusion for any
integer $r$, and all sufficiently large $m$ (depending on $r$). Consider
the pairs $(r,m)$ for which the conclusion of the Proposition
holds. Given $r$, we have seen already that there is a positive integer
$r_0>r$ so that (i) the conclusion holds for $(r_0,m)$ for all $m\geq
m_0$, say, and (ii)
so that $\Oh_{\tilde{X}}(-(r_0-r)E)$ is very ample for $\pi$. Then
choose $m_1$ so that  
$\pi_*\Oh_{\tilde{X}}(-(r_0-r)E)\tensor_{\Oh{X}}{\mathcal L}^{\tensor
m}$ is globally generated for all $m\geq m_1$. Then
$\Oh_{\tilde{X}}(-(r_0-r)E)\tensor{\mathcal M}^{\tensor m}$ is globally
generated, for all $m\geq m_1$. Hence the conclusion of the proposition
holds for $(r,m)$ with $m\geq m_0+m_1$.
\end{remark}

\begin{prop}\label{cokernel}
For $\tilde{X}$, $\tilde{Y}$ as above, the condition $\ALeff(\tilde
X,\tilde Y)$ holds.
\end{prop}
\begin{proof}
By Proposition~\ref{prop1}, for any invertible formal sheaf $\scrF$ on
$\frakX$, one has a map of formal locally free sheaves
\begin{equation}\label{eqsurj1}
({\hat\scrM}^{\tensor{-M}})^{\oplus s}(-F_1)
 \to \scrF \to 0 \hspace{5mm}  
\end{equation}
for some $M>>0$, $s>0$, with cokernel supported in
$\tilde{Y}\setminus E$, where $F_1$ is an effective divisor on
$\tilde{X}$ with exceptional support.

Similarly, for the dual formal line bundle $\scrF^{\vee}$, we have a
map, sujective outside the exceptional locus,
$$(\hat\scrM^{\tensor{-N}})^{\oplus t}(-F_2)
\to \scrF^{\vee} \to 0 \hspace{5mm}$$ 
for some $N>>0$, $t>0$. Dualizing this we have an injection 
\begin{equation}\label{eqsurj2}
\scrF \to (\hat\scrM^{\tensor N})^{\oplus
t}(F_2)
\end{equation}
which is a split inclusion on stalks at any point of $\tilde{Y}\setminus
E$.

Composing the maps in (\ref{eqsurj1}) and (\ref{eqsurj2}), we have a map
between formal locally free sheaves
\begin{equation}\label{eq3}
(\hat\scrM^{\tensor{-M}})^{\oplus s}(-F_1)
\by{\hat\tau}
(\hat\scrM^{\tensor{N}})^{\oplus t}(F_2)
\end{equation}
such that $\im{\hat\tau} \into \scrF$ with cokernel supported in
$\tilde{Y}\setminus E$. 

The map $\hat{\tau}$ may be described by an $s\times t$ matrix of elements
of 
\[H^0(\frakX,\widehat{\scrM}^{\tensor N+M}(F_1+F_2)).\]
By the condition $\Lef(\tilde{X},\tilde{Y})$ 
(or rather Corollary~\ref{formalsects}), 
\[H^0(\frakX,\widehat{\scrM}^{\tensor
N+M}(F_1+F_2))\cong H^0(\tilde X,\scrM^{\tensor N+M}),\]
so that the map $\hat\tau$ is the formal completion of a map of locally
free sheaves on $\tilde{X}$
$$(\scrM^{\tensor{-M}})^{\oplus s}(-F_1) \by{\tau}
(\scrM^{\tensor{N}})^{\oplus t}(F_2).$$ 
Thus we have $\widehat{\im(\tau)} \into \scrF$, with cokernel
supported on $\tilde{Y}\setminus E$. 

Now $\im(\tau)$ is a coherent sheaf on $\tilde{X}$, such that for any
point $y\in \tilde{Y}\setminus E$, the stalk at $y$ satisfies
\[\im(\tau)_y\tensor_{\Oh_{\tilde{X},y}}\widehat{\Oh_{\tilde{X},y}}\cong
\scrF_y\cong\Oh_{\frakX,y}\cong\widehat{\Oh_{\tilde{X},y}},\]
(where the completions are with respect to the ideal defining $\tilde Y$).
Hence $\im(\tau)$ is a coherent sheaf of rank 1, which is invertible on
$\tilde{X}$ at all points in $\tilde{Y}\setminus E$. Since $\tilde{X}$ 
is non-singular, the double dual of $\im(\tau)$ is an invertible sheaf
$\widetilde{\scrF}$ on $\tilde{X}$, such that 
\[\widetilde{\scrF}\tensor{\Oh_{\tilde Y}}\mid_{\tilde{Y}\setminus
E}\cong
{\scrF_1}\mid_{\tilde{Y}\setminus E}.\]
Thus we have $\ALeff(\tilde{X},\tilde {Y})$.
\end{proof}

\begin{remark} The above argument, applied to an arbitrary formal locally 
free sheaf $\scrF$, implies the existence of a coherent, reflexive sheaf
$\mathcal F$ together with an injective map $\widehat{\mathcal F}\to
\scrF$ which restricts to an isomorphism on $\tilde{Y}\setminus E$. We do 
not know if $\mathcal F$ can be chosen to be locally free in a
neighbourhood of $\tilde{Y}$; perhaps the ``natural'' extension to
our situation of the Grothendieck ``Leff'' condition is for this property
to hold.
\end{remark}

\begin{cor}\label{cok}
 If $\dim \tilde X\geq 4$, the cokernel of the restriction map
$\Pic(\tilde X) \to \Pic(\tilde Y)$ is generated by exceptional divisors 
of $\tilde Y$ which map to points in $Y$.
\end{cor}
\begin{proof}
By Lemma \ref{exp} one has $\Pic(\frakX) \isom \Pic(\tilde
Y)$. Let $\mathcal N$  be a line bundle  on $\tilde Y$, and $\mathfrak N$
its unique lift to a formal line bundle on $\frakX$. Proposition
\ref{cokernel} implies that there exists an invertible sheaf $\mathcal G$
on $\tilde X$ such that $\hat{\mathcal G} \isom \mathfrak N$ on
$\tilde{Y}\setminus E$.
Thus ${\mathcal G}\mid_{\tilde{Y}}\tensor{ \mathcal N}^{\vee}$ is a
line 
bundle on $\tilde{Y}$ which has trivial restriction to $\tilde{Y}\setminus
E$, and is thus the line bundle associated to a divisor on $\tilde{Y}$
with exceptional support. It remains to show that, upto tensoring with a
line bundle restricted from $\tilde{X}$, it corresponds to a sum of
exceptional divisors for $\tilde Y\to Y$ with 0-dimensional image in $Y$. 

Let $E_1,\ldots,E_r$ be the irreducible exceptional divisors of
$\tilde{X}\to X$, indexed so that for some $0\leq
s\leq t\leq r$, we have
\begin{enumerate} 
\item[(i)] $E_1,\ldots,E_s$ are the irreducible divisors in
$\tilde{X}$ with 0-dimensional image in $X$
\item[(ii)] $E_{s+1},\ldots,E_t$ are the irreducible divisors with 
1-dimensional image in $X$
\item[(iii)] $E_{t+1},\ldots,E_r$ are the irreducible exceptional divisors
for $\tilde{X}\to X$ whose images in $X$ have dimension $\geq 2$. 
\end{enumerate}

Since $\tilde Y\in |\VV|$ is a general member, we have 
\begin{enumerate}
\item[(a)] $\tilde{Y}\cap E_i=\emptyset$ for $1\leq i\leq s$
\item[(b)] for each $s+1\leq i\leq t$, let $\tilde{Y}\cap
E_i=\cup_{j=1}^{s_i}F_{ij}$ be the irreducible (equivalently
connected) components of the intersection; then $\tilde{Y}\cap E_i$ is
non-singular, reduced, and has no common irreducible component with
$\tilde{Y}\cap E_{i'}$ for any $i'\neq i$
\item[(c)] for $t+1\leq i\leq r$,
$F_i=\tilde Y\cap E_i$ is reduced and irreducible.
\end{enumerate}
Here, (a) is clear, and (c) follows from Bertini's theorem. For (b), note
that the linear system $|\VV|$ restricts to a base-point free linear
system $|\VV_i|$ on $E_i$, and the Stein factorization of the
corresponding morphism has the form $\pi_i:E_i\to C_i$ for some
irreducible curve $C_i$, such that $|\VV_i|$ is the pull-back of a linear
system from $C_i$. Hence the general member of $|\VV_i|$ is a disjoint
union of a finite set of fibers of $E_i\to C_i$, over points of $C_i$ for
which $E_i\to C_i$ is smooth. 

Thus, the line bundle determined by an irreducible exceptional divisor $F$
on $\tilde{Y}$ lies in the image of the $\Pic(\tilde{X})\to\Pic(\tilde 
Y)$, except possibly when $F$ is one of the $F_{ij}$ in (b) above, and the
image of each such $F_{ij}$ in $Y$ is a point.
\end{proof}
\begin{lemma}\label{independence}
Assume $\dim \tilde X\geq 3$, and $\tilde Y\in|\VV|$ is general. Let $F$
be a divisor on $\tilde{Y}$ with exceptional support, such that 
$\Oh_{\tilde{Y}}(F)$ is the restriction of a line bundle from $\tilde{X}$. 
Then there is a divisor $\tilde F$ on $\tilde{X}$ with exceptional support 
such that $\Oh_{\tilde X}(\tilde{F})\mid_{\tilde{Y}}=\Oh_{\tilde{Y}}(F)$.
\end{lemma}
\begin{proof}
From the description of the irreducible exceptional divisors of
$\tilde{Y}\to Y$ given in (a), (b), (c) of the proof of
Corollary~\ref{cok} above, we see that it suffices to assume (in the
notation of (b)) that
\begin{equation}\label{eqexc}
F=\sum_{i=s+1}^t\sum_jn_{ij}F_{ij},
\end{equation}
for some integers $n_{ij}$, $s+1\leq i\leq t$, $1\leq j\leq s_i$. Here,
the divisors
\[F_i:=\sum_{j=1}^{s_i}F_{ij}\]
for each $s+1\leq i\leq t$ satisfy $\Oh_{\tilde{X}}(E_i)\mid_{\tilde
Y}\cong \Oh_{\tilde Y}(F_i)$. 

So it suffices to show: if $F$ in (\ref{eqexc}) is such that
$\Oh_{\tilde{Y}}(F)$ is the restriction of a line bundle from $\tilde{X}$,
then $n_{ij}$ is independent of $j$, for each $s+1\leq i\leq t$. This is
done using a suitable computation with intersection numbers.

Let $T\subset \tilde{X}$ be a general complete intersection of dimension
3, in some projective embedding of $\tilde{X}$. Then by Bertini's theorem,
we may assume that 
\begin{enumerate}
\item[(i)] $T$ is irreducible and nonsingular, and the scheme
theoretic intersection $T\cap E_i$ is a reduced, irreducible surface, for
each $1\leq i\leq r$ 
\item[(ii)] the scheme theoretic intersection $Z:=T\cap\tilde{Y}$ is a
(reduced, irreducible) nonsingular surface in $T$ 
\item[(iii)] $Z\cap E_i=\emptyset$ for $1\leq i\leq s$, and $Z\cap E_i$ is
a reduced, irreducible curve for each $t+1\leq i\leq r$
\item[(iv)] $Z\cap F_{ij}$ is a reduced, irreducible curve in $Z$, for
each $s+1\leq i\leq t$, for all $j$.
\end{enumerate}

If $T_0$ is the image of $T$ in $X$, let $\overline{T}\to T_0$ be the
normalization. Then $\pi_T:T\to \overline{T}$ is a resolution of 
singularities, such that $T\cap E_i$, $1\leq i\leq r$ are the irreducible
exceptional divisors of $\pi_T$. The surface $\pi_T(Z)=\overline{Z}\subset
\overline{T}$ is a normal Cartier divisor in $\overline{T}$ (it is a
general member of the ample, base-point free linear system on
$\overline{T}$ determined by $|\VV|$). Let $\pi_Z:Z\to\overline{Z}$ be the
restriction of $\pi_T$; then $\pi_Z$ is a resolution of singularities of a
normal, projective surface, with irreducible exceptional curves 
$Z\cap E_i$, $t+1\leq i\leq r$ and $Z\cap F_{ij}$, $s+1\leq i\leq t$, 
$1\leq j\leq s_i$. 

Since the linear equivalence class of $F=\sum_{i,j}n_{ij}F_{ij}$ is
assumed to lie in the image of the restriction
$\Pic\tilde{X}\to\Pic\tilde{Y}$, we have that $\Oh_Z(\sum_{i,j}n_{ij}Z\cap
F_{ij})$ is in the image of $\Pic T\to \Pic Z$. For each $i$, the $F_{ij}$
are irreducible components of general fibers of $E_i\to C_i$. Hence
$F_{ij}\cap Z$ are irreducible components of general fibers of $T\cap
E_i\to C_i$, and are thus {\em algebraically equivalent} as 1-cycles on
the smooth projective 3-fold $T$. Hence, for any divisor $D$ on $T$, the
intersection number $(D\cdot (F_{ij}\cap T))_T$ is independent of $j$, for
each $i$. Since $F_{ij}\cap Z$ is a Cartier divisor in $Z$,  
if $D_Z$ is any divisor on $Z$ representing $\Oh_T(D)\mid_Z$,
the projection formula gives an equation between
intersection numbers
\begin{equation}\label{intnum}
(D_Z\cdot (F_{ij}\cap Z))_Z=(D\cdot (F_{ij}\cap T))_T
\end{equation}
computed on the surface $Z$ and the 3-fold $T$, respectively.

Apply this to our divisor $F\mid_Z=\sum_{i,j}n_{ij}F_{ij}\cap Z$,
which is assumed to be of the form $D_Z$ for some divisor $D$ on $T$. We
get that
\begin{equation}\label{intnum2}
(\sum_{ij}n_{ij}(F_{ij}\cap Z)\cdot (F_{i'j'}\cap Z))_Z
\end{equation}
is independent of $j'$, for each $s+1\leq i'\leq t$.

Since $F_{ij}\cap Z$ are irreducible exceptional divisors for
$\pi_Z:Z\to\bar{Z}$, which is a resolution of singularities of a
normal surface, the intersection matrix 
\[(F_{ij}\cap Z,F_{i'j'}\cap Z)_Z\]
is negative definite. Regard (\ref{intnum2}) as a system of linear
equations satisfied by the $n_{ij}$, with coefficients given by
intersection numbers.  The solutions of the system
(\ref{intnum2}) for the ``unknowns'' $n_{ij}$ correspond to elements in
the $\Z$-span of the $F_{ij}$, which are in the orthogonal complement of
the span of all the differences $(F_{ij_1}\cap Z)-(F_{ij_2}\cap Z)$, for
all $1\leq j_1<j_2\leq s_i$, and all $s+1\leq i\leq t$. The span of 
these differences clearly has co-rank $t-s$, so the orthogonal
complement has rank $t-s$. We have $t-s$ elements 
$\sum_{j=1}^{s_i}(F_{ij}\cap Z)$ which lie in the orthogonal complement,
which are clearly independent, so must span the orthogonal complement  
after tensoring with $\mathbb Q$. This implies that 
$F_Z=\sum_{i,j}n_{ij}(F_{ij}\cap Z)$ must be a rational linear combination
of the divisors $\sum_{j=1}^{s_i}(F_{ij}\cap Z)$, and so $n_{ij}$ must be
independent of $j$, for each $i$, as desired.
\end{proof}

Assume now that $\dim X\geq 4$. The conditions $\Lef$ and $\ALeff$ imply
that we have the following diagram, with exact rows and columns:
\[
\begin{diagram}{ccccccccccc}
& &  && 0 && 0 &&  & & \\
&  &  & &  \downarrow & & \downarrow & &  & & \\
0 & \to & I_{X}^{'} & \to &\Z[E_{X}] & \to & \Z[E_{Y}] & \to & I_{Y}^{'}
&
\to & 0 \\
& & \downarrow &  & \downarrow & & \downarrow & & \downarrow & & \\
0 & \to & I_{X} & \to & \Pic(\tilde X) & \to & \Pic(\tilde Y) & \to & I_{Y} 
& \to & 0 \\
&  & & &  \downarrow & & \downarrow & & & & \\
& & &  & \Cl(X) & \to & \Cl(Y)  &  & & & \\
&  & & &  \downarrow & & \downarrow & & & & \\
& & &  & 0 & & 0  & & & & \\
\end{diagram}
\]
Here $\Z[E_{X}]$ and $\Z[E_{Y}]$ are the subgroups in the respective
Picard groups freely generated by the irreducible exceptional divisors in 
$\tilde{X}$ and $\tilde{Y}$, and $I_X$, $I'_X$, $I_Y$, $I'_Y$ are defined
by the exactness of the rows. Clearly $I'_{X}$, $I_{Y}^{'}$ are 
generated by the irreducible exceptional divisors in $\tilde{X}$ and
$\tilde{Y}$, respectively, which have 0-dimensional image
under $\pi$.
  
It is clear that $I_{X}\isom I_{X}^{'}$: it is an
injection since $\Z[E_{X}] \into \Pic(\tilde X)$ is so. That it is a
surjection follows from Corollary~\ref{kernel}. Also Corollary \ref{cok}
shows that that the map $I_{Y}^{'} \to I_{Y}$ is surjective, while it is
also injective, by Lemma~\ref{independence}. Hence, from a diagram chase,
we see that $\Cl(X)\to \Cl(Y)$ is an isomorphism, completing the proof of 
Theorem~\ref{thm1'} when $\dim \tilde{X}\ge 4$. 

By a similar argument, Corollary~\ref{kernel} and Lemma~\ref{independence}
imply that, if $\dim \tilde X=3$, then  $\Cl(X)\to \Cl(Y)$ is
{\em injective}. The finite generation of the cokernel results from the
fact that $\Pic^0\tilde{X}\to\Pic^0\tilde Y$ is an isogeny (and hence an
isomorphism), since the map on tangent spaces $H^1(\tilde X,\Oh_{\tilde
X})\to H^1(\tilde Y,\Oh_{\tilde Y})$ is an isomorphism, from
Theorem~\ref{GR}. This completes the proof of Theorem~\ref{thm1'} in case
$\dim \tilde X=3$.

For possible use elsewhere, we make explicit the following result, more
or less implicit above. We thank the referee for some illuminating remarks
about formal Cartier divisors. 
\begin{thm}\label{formal} 
Let $\tilde{X}$ be as in lemma~\ref{redn}, with $\dim \tilde{X}= 3$, and
let $E$ be an effective divisor on $\tilde{X}$ with exceptional support
such that $-E$ is $\pi$-ample. Let $\tilde{Y}\subset \tilde{X}$ be a
general member of the linear system $|\VV|$, and let $\frakX$ denote the
formal completion of $\tilde{X}$ along $\tilde{Y}$.  Then the map
\[\rho_{\frakX}:\Pic\tilde{X}\to\Pic \frakX\]
has the following properties.
\begin{enumerate}
\item [(i)] The kernel of $\rho_{\frakX}$ is freely generated by the
classes of irreducible $\pi$-exceptional divisors with 0-dimensional image
under $\pi$.
\item[(ii)] The cokernel of $\rho_{\frakX}$ is generated by the classes of
exceptional divisors $F$ on $\tilde{Y}$ such that $\dim \pi({\rm
supp}\,F)=0$ (in particular, the corresponding line bundles on $\tilde{Y}$
do extend to formal line  bundles).
\item[(iii)] With the notation introduced above (proof of
Corollary~\ref{cok}),  let $A$ denote the quotient of the free abelian
group on $F_{ij}$, $s+1\leq i\leq t$, $1\leq j\leq s_i$, by the subgroup   
generated by $\sum_{j=1}^{s_i}F_{ij}$, for $s+1\leq i\leq t$. Then there
is a natural isomorphism $\coker(\rho_{\frakX})\cong A$.
\item[(iv)] Let $\pi:\tilde{X}\to
X$, $\pi:\tilde{Y}\to Y$ be the Stein factorizations. There is a natural
isomorphism  
\[\coker({\rm Cl}\,(X)\to{\rm Cl}\,(Y))\cong \coker ({\rm
Pic}\,(\frakX)\to {\rm Pic}\,(\tilde{Y})).\]
\end{enumerate}
\end{thm}
Most of these conclusions have already been obtained in the course of the
above proof. The only remaining assertion to prove is that all the
line bundles $\Oh_{\tilde{Y}}(F_{ij})$ do extend to formal line bundles on
$\frakX$. We thank the referee for suggesting a proof.

In fact, for each $i$ as in (iii), the divisor $E_i$ defines a 
formal Cartier divisor $\hat{E_i}$ on $\frakX$ with support $E_i\cap
\tilde{Y}$, which is a disjoint union of the closed sets $F_{ij}$, $1\leq
j\leq s_i$. Hence each connected component of the support defines a formal
Cartier divisor $\hat{F_{ij}}$ on $\frakX$ (one may first check a similar
assertion for Cartier divisors on each of the schemes $\frakX_m$, all of
which have reduced scheme $\tilde{Y}$, for example). The associated formal
line bundles are the desired extensions. 
\end{subsection}

\end{section}

\begin{section}{Application to $1$-motives}
In \cite{D}, \S10, Deligne defined 1-motives over $k$ as complexes $[L\to
G]$, where $L$ is a lattice (free abelian group of finite rank with a
continuous action of the absolute Galois group of $k$), and $G$ a
semi-abelian $k$-variety. This gives an algebraic way of ``defining''
certain
(co)homology groups of a variety, in a manner analogous to the way in
which the Jacobian of a non-singular projective curve ``defines'' its
first (co)homology group algebraically. Over $\C$, 1-motives have a
transcendental description using certain special types of mixed Hodge
structures, and there is an equivalence of categories between 1-motives
over $\C$ and the full subcategory of these special types of mixed Hodge
structures. In particular, there is an underlying philosophy (``Deligne's
Conjecture'', some aspects of which have been proved in \cite{BRS},
\cite{R}) that, if one can construct a 1-motive transcendentally using
some ``part'' of the mixed Hodge structure of an algebraic variety, then
there must be an algebraic construction of that 1-motive as well, valid
over more general ground fields. Further, if some operation between
1-motives can be constructed transcendentally, there must be an algebraic
construction of it as well, and properties of such an operation
(e.g. injectivity, isomorphism) should have algebraic proofs.  

In \cite{BS}, a 1-motive $\Alb^+(X)$, the {\em cohomological Albanese
1-motive}, has been associated to any variety $X$ over a field $k$ of
characteristic 0. If $X$ is proper, this is a semi-abelian variety over
$k$, and if $X$ is also non-singular, it coincides with the
``classical'' Albanese variety. If $k=\C$, $\Alb^+(X)$ can be
constructed analytically, using the mixed Hodge structure on
$H^{2n-1}(X,\Z(n))$, where $n=\dim X$, in a manner generalizing the
analytic construction of the Albanese variety of a non-singular
proper complex variety. For a proper, possibly singular complex variety
$X$, one has a formula
\[\Alb^+(X)(\C)=\Ext^1_{\bf MHS}(\Z,H^{2n-1}(X,\Z(n)))\]
where the right side is the group of extensions in the (abelian) category
of mixed Hodge structures. 

If $X$ is projective, $Y\subset X$ is a (reduced, effective) Cartier
divisor, then there is a {\em Gysin map} $\Alb^+(Y)\to \Alb^+(X)$,
constructed algebraically in \cite{BS}, and which in case $k=\C$ corresponds 
to the Gysin map (modulo torsion) $H^{2n-3}(Y,\Z(n-1))\to
H^{2n-1}(X,\Z(n))$ in topology (which is a morphism of mixed Hodge
structures). 

In case $X$ is projective over $\C$, and $Y$ is a general hyperplane
section (here, ``general'' means ``in a Zariski open set of the
parameter variety''), then it is shown in \cite{BiS} that $\Alb^+(Y)\to
\Alb^+(X)$ is an isomorphism when $\dim X=n\geq 3$; this is an important
step in the proof of the Roitman theorem for singular projective
complex varieties (the main result of \cite{BiS}). In case $X$ (and hence
$Y$) is non-singular, this is a particular case of the Lefschetz
hyperplane theorem.  The proof of this isomorphism in \cite{BiS} is
transcendental, ultimately relying on the local structure (in the
Euclidean topology) of a morphism of complex varieties, which is given by
the theory of Whitney stratifications (see \cite{Ve}, or \cite{GM}, for
example). 

This suggests that, if $X$ is projective over a field $k$ of
characteristic 0, of dimension $\geq 3$, and $Y\subset X$ is a general
hyperplane section, then the Gysin map $\Alb^+(Y)\to\Alb^+(X)$ is an
isomorphism; further, the ``philosophy of 1-motives'' suggests that there
is a purely algebraic proof of this fact. 

The validity of the isomorphism over an arbitrary $k$ of characteristic 0
can be deduced from the case $k=\C$. An algebraic proof, on the other
hand, can be obtained as follows. 

It is easy to see that the general case follows from the case when $k$ is
algebraically closed, so we assume this holds. Next, the category of
1-motives admits a notion of Cartier duality, which is an
auto-antiequivalence of the category. So it suffices to show that the
Cartier dual to $\Alb^+(Y)\to \Alb^+(X)$ is an isomorphism. 

It is shown in \cite{BS} that the Cartier dual of $\Alb^+(X)$ is another
1-motive, explicitly described as follows (implicitly, this gives an
algebraic description of $\Alb^+(X)$). 

Let $\pi:X'\to X$ be the normalization map, $\Cl(X')$ the divisor class
group of the normal projective variety $X'$, and $\Cl^0(X')$ the largest
divisible subgroup. Then $\Cl^0(X')$ is naturally identified with (the
$k$-points of) an abelian variety (which can be identified with the Picard
variety of any resolution of singularities).  Let $L_X$ denote the group
of all Weil divisors $D$ on $X'$ such that 
\begin{enumerate}
\item[(i)] $\pi_*(D)=0$ as a cycle on $X$ 
\item[(ii)] $[D]\in \Cl(X)$ lies in the subgroup $\Cl^0(X')$. 
\end{enumerate}
If ${\rm Div}(X'/X)$ denotes the group of Weil divisors $D$ on $X'$ with
$\pi_*(D)=0$ as a cycle on $X$, then 
\begin{equation}\label{eq1}
L_X=\ker \left(\Div(X'/X)\to \Cl(X')/\Cl^0(X')\right).
\end{equation}
Note that ${\rm Div}(X'/X)$ is a subgroup of the group
of Weil divisors on $X'$ which are supported on $\pi^{-1}(X_{sing})$; in
particular, $\Div(X'/X)$ is free abelian of finite rank.
Thus $L_X$ is a free abelian group of finite rank, and the
obvious homomorphism $L_X\to \Cl^0(X')$ defines a 1-motive; this is the
Cartier dual to $\Alb^+(X)$. 

We are given that $Y$ is a general hyperplane section of $X$ in a certain
projective embedding; the pull-back of the corresponding very ample linear
system to $X'$ gives an ample, base-point free linear system on
$X'$. Hence $Y'=Y\times_XX'$ is a general member of this linear system on 
$X'$, and is thus normal, by Bertini's theorem, so that $Y'\to Y$
is the normalization of $Y$. There is also an associated restriction map
$\Cl(X')\to \Cl(Y')$, such that by Theorem~\ref{thm1} above, it is an
isomorphism when $n=\dim X\geq 4$, and is an injection with finitely
generated cokernel if $n=3$. Hence it induces an isomorphism $\Cl^0(X')\to
\Cl^0(Y')$ between abelian varieties, and an inclusion on quotients
$\Cl(X')/\Cl^0(X')\to \Cl(Y')/\Cl^0(Y')$, if $\dim X\geq 3$. 

Now the Gysin map $\Alb^+(Y)\to\Alb^+(X)$ is Cartier dual to a map of
1-motives (i.e., to a map between 2-term complexes)
\[\left[L_X\to \Cl^0(X')\right]\;\;\longrightarrow\;\;\left[L_Y\to
\Cl^0(Y')\right]\]
where $\Cl^0(X')\to\Cl^0(Y')$ is the above isomorphism, induced by the
restriction homomorphism $\Cl(X')\to \Cl(Y')$. It remains to see that, for
general $Y$, this map is an isomorphism of 1-motives, i.e., the 
map $L_X\to L_Y$ is also an isomorphism. 

Since $Y'=Y\times_XX'$, a functorial property of the refined Gysin
homomorphism defined in \cite{Fulton} implies that if $D$ is a Weil
divisor on $X'$ with $\pi_*(D)=0$ as a cycle on $X$, then the cycle
$[D\cap Y']$ has the property that $\pi_*[D\cap Y']=0$ as a cycle on
$Y$. This gives us a map $L_X\to L_Y$. It is shown in \cite{BS} that this
is the map corresponding to the Cartier dual of the Gysin map. 

Since $Y$ is a general member of a very ample linear system on $X$,
where $\dim X\geq 3$, it is clear (using Bertini's theorem) that if $D$ is
an irreducible Weil divisor in $Y'$ lying over the singular locus of $Y$,
then there is a unique irreducible Weil divisor $D_1$ in $X'$, lying over
the singular locus of $X$, such that $D=D_1\cap Y'$ as divisors. This is
because $Y_{sing}=X_{sing}\cap Y$, giving a bijection between the
codimension 1 irreducible components of $X_{sing}$ and $Y_{sing}$, which
also gives a bijection between the codimension 1 irreducible components of 
$\pi^{-1}(X_{sing})\subset X'$ and $\pi^{-1}(Y_{sing})\subset Y'$. 
Hence, for general $Y$ as above, we have that ${\rm Div}(X'/X)\to {\rm
Div}(Y'/Y)$ is an isomorphism. The formula (\ref{eq1}) applied to $X$ and
to $Y$, together with the fact (from Theorem~\ref{thm1}) that
$\Cl(X')/\Cl^0(X')\to\Cl(Y')/\Cl^0(Y')$ is injective, implies that $L_X\to
L_Y$ is an isomorphism, as was to be shown.  

\end{section}

\section{Some refinements, and statements in any characteristic}

In this brief section, we make connections with the classical theory of
Weil divisors, in the style of Weil \cite{Weil} and Lang, as exposed in 
Lang's book \cite{Lang}, following comments of the referee. This gives
another perspective to the above results, and yields also some statements 
in arbitrary characteristic. 

First, in Theorem~\ref{thm1}, one can improve the statement in the
following ways.
\begin{enumerate}
\item[(i)] There is a dense Zariski open set $\Omega\subset
|\VV|$ so that, if $K$ is any algebraically closed extension field of $k$,
$X_K=X\times_kK$, and $Y_K\subset X_K$ is a memeber of the base-changed
linear system $|\VV\tensor_kK|$ on $X_K$, corresponding to a $K$-point of 
$\Omega$, then the theorem holds for the pair $(X_K,Y_K)$. As stated, 
Theorem~\ref{thm1} does yield such an open subset of $|\VV\tensor_kK|$,
but in fact it may be taken to be $\Omega_K$. It is not difficult to
modify the proof given above to yield this conclusion as well.
\item[(ii)] When $\dim X=3$, the cokernel of the (injective) map on class
groups is in fact torsion-free. This follows from the proof given, since
from Theorem~\ref{formal}, it boils down to the assertions that the
cokernel of $\Pic \frakX\to \Pic \tilde{Y}$ is torsion-free. If ${\mathcal
K}=\ker({\mathcal O}_{\frakX}^*\to{\mathcal O}_{\tilde{Y}}^*)$, then
$\mathcal K$ is a sheaf of $\mathbb Q$-vector spaces on the topological
space $\tilde{Y}$ (as may be immediately verified on suitable affine
open subsets), and we have an exact sequence
\[\Pic\frakX \to \Pic \tilde{Y}\to H^2(\tilde{Y},{\mathcal K})\]
where the last term is a $\mathbb Q$-vector space.
\end{enumerate}
A different way of seeing (ii) is by a transcendental argument,
using a suitable Lefschetz theorem, as in the Appendix: in
Theorem~\ref{smt}, for $i=\hat{n}$, the cokernel is torsion-free and the
conclusions of that theorem also hold for cohomology with 
${\mathbb Z}/n{\mathbb Z}$ coefficients, for any $n>0$. 

The referee has also pointed out that for a geometrically integral  
projective variety $X$ over a field $k$, which is smooth in codimension 1,  
one can associate to it an abelian $k$-variety $\Pic_W(X)$ (the ``Picard
variety in the sense of Weil''), such that when $k=\bar{k}$, the group of
$k$-rational points $\Pic_W(X)(k)$ coincides with the group ${\rm
Cl}^0(X)$, the maximal divisible subgroup of the group ${\rm Cl}\,(X)$ of
Weil divisors modulo linear equivalence. Further, the Weil-Neron-Severi
group $NS_W(X)={\rm Cl}\,(X)/{\rm Cl}^0\,(X)$ is finitely generated.  
With this notation, we sketch the referee's argument to prove the
following result. 
\begin{thm}\label{referee} 
Let $X$ be an irreducible projective $k$-variety of dimension $d\geq 3$,
which is regular in codimension 1, where $k$ is an algebraically closed
field of any characteristic. Let $f:X\to \P^N_k$ be an embedding. Let
$\bar{K}$ be the algebraic closure of the function field
$K=k(\hat{\P}^N_k)$ of the dual projective space, and
let $Y_K\subset X_K$ be the generic hyperplane section. Then
\begin{enumerate}
\item[(i)]
$\Pic_W(X)_K\to \Pic_W(Y_K)$ is an isomorphism 
\item[(ii)] the composition 
\[ NS_W(X)\to NS_W(X_{\bar{K}})\to NS_W(Y_{\bar{K}})\] 
is injective.
\end{enumerate}
\end{thm}
The isomorphism between Weil-Picard varieties is a consequence of
\cite{Lang}~VIII, Theorem~4. The injectivity on Weil-Neron-Severi groups
is reduced to a result of Weil~\cite{Weil}. We must show that, if
$\mathcal 
L$ is a line bundle on $X$, whose pull-back to $(Y_K)_{\rm reg}$, the
smooth locus of $Y_K$, is algebraically equivalent to 0, then $\mathcal L$
is algebraically equivalent to 0 on $X$. From \cite{Lang}~VI, Theorem~1,
the pullback of $\mathcal L$ to $(Y_K)_{\rm reg}$ determines a
$K$-rational point of $\Pic_W(Y_K)=\Pic_W(X)_K$. Since $K$ is a pure
transcendental extension of $k$, this must determine a $k$-rational point
of $\Pic_W(X)$. Hence, changing $\mathcal L$ by  the class of some point
of $\Pic_W(X)$, we may assume $\mathcal L$ has trivial pull-back to
$(Y_K)_{\rm reg}$. Now Theorem~2 of Weil~\cite{Weil} implies that
$\mathcal L$ is itself trivial on $X$.

\begin{section}{Appendix: Grothendieck-Lefschetz theorem for complex
 projective varieties}

In this appendix, we shall sketch the proof of the following theorem
using results from stratified Morse theory, as explained to us by
Najmuddin Fakhruddin.
\begin{thm}\label{GL2}
Let $X$ be a smooth projective variety of dimension at least $4$
defined over the field of complex numbers, $\bbC$. Let $\scrL$ be a
big line bundle over $X$ generated by global sections. If $Y$ denotes
a general member of the linear system $\vert\scrL\vert$, then one has
an exact sequence
$$ 0 \to K \to \Pic(X) \to \Pic(Y) \to Q \to 0 $$ where $K$ is the
(free) subgroup generated by divisors in $X$ which map to points under
the generically finite map $X \by{\pi} \bbP(\HH^{0}(X,\scrL))$ and $Q$
is the group generated by the irreducible components of the
restriction of divisors in $X$ which map to curves under $\pi$.
\end{thm}
Theorem~\ref{GL2} is an immediate consequence of
Corollaries~\ref{Pic^0cor} and \ref{NScor} below.
 
In what follows, all cohomologies that we consider are singular
cohomology of the underlying analytic space(s), with $\bbZ$-coefficients.
Recall that, for any $\mathbb C$-variety, these cohomology groups support
mixed Hodge structures, which are functorial for morphisms
between varieties (see \cite{D}). The proofs of Corollaries~\ref{Pic^0cor}
and \ref{NScor} are reduced to assertions about the homomorphisms
$\HH^i(X)\to \HH^i(Y)$ for $i=1,2$, using the following standard lemma.
\begin{lemma}\label{red-to-top}
Let $W$ be a smooth proper $\mathbb C$-variety. Then there are
isomorphisms, functorial for morphisms of $\mathbb C$-varieties,
\[\Pic^0(W)\cong\frac{H^1(W)\tensor {\mathbb C}}
{F^1H^1(W)\tensor{\mathbb C}+H^1(W)},\]
\[NS(W)=\ker\left(H^2(W)\to\frac{H^2(W)\tensor{\mathbb
C}}{F^1H^2(W)\tensor{\mathbb C}}\right).\]
\end{lemma}
\begin{proof} From Serre's GAGA, it follows that
$\Pic(W)\cong\Pic(W_{an})$, where the latter denotes the group of
isomorphism classes of analytic line bundles. Using the exponential sheaf
sequence, and the Hodge decomposition, we obtain the above isomorphisms in
a standard way, where $\Pic^0(W)=\ker(\Pic(W)\to H^2(W))$ is the maximal
divisible subgroup, and the Neron-Severi group ${\rm NS}\,(W)={\rm
image}\,(\Pic(W)\to \HH^2(W))$ is finitely generated.
\end{proof}

We now state the following consequence of the Relative
Lefschetz theorem with Large fibres (see \cite{GM}, page 195). 

\begin{thm}\label{smt}
Let $W$ be a $n$-dimensional nonsingular connected algebraic
variety. Let $\pi:W \to \bbP^{N}$ be a morphism and let $H \subset
\bbP^{N}$ be a general linear subspace of codimension $c$. Define
$\phi(k)$ to be the dimension of the set of points $z \in \bbP^{N}$
such that the fibre $\pi^{-1}(z)$ has dimension $k$. (If this set is
empty define $\phi(k)=-\infty.$) Then the homomorphism induced by
restriction,
$$ \HH^{i}(W,\bbZ) \to \HH^{i}(\pi^{-1}(H), \bbZ)$$ is an isomorphism for
$i<\hat{n}$ and is an injection for $i=\hat{n}$, where
$$\hat{n}=~n-\underset{k}{sup}~(2k -(n-\phi(k)) + inf(\phi(k),c-1)) -1 $$
\end{thm}

In the situation of theorem \ref{GL2}, we first take $W=X$, $\pi$ a
generically finite map and $H$ a general hyperplane. Then one can
easily check that $\hat{n}\ge 1$ in this case. Let $X'$ be the (open)
subvariety of $X$ defined by removing all divisors which map to
points under $\pi$, and $X''$ be the subvariety obtained by further
removing divisors which map to curves under $\pi$. In these two
cases, for the restriction of $\pi$ to $X'$ and $X''$, one can check
that $\hat{n}\ge 2$ and $\ge 3$ respectively.

Let $Y=\pi^{-1}(H)$, and let $Y'$ and $Y''$ be defined similarly
in $X'$ and $X''$ respectively. Since a general hyperplane section in
$\bbP^{N}$ misses points, one notes immediately that $Y'=Y$.

\begin{lemma}\label{tf} If $V\subset W$ is a dense Zariski open subset of
a non-singular proper variety $W$, then 
\begin{enumerate}
\item[(i)] $\HH^1(W,V)=0$, and $\HH^{2}(W,V)$ is a free abelian group,
with a basis given by the irreducible divisors supported on $W\setminus V$
(in particular, it is pure of weight 2)
\item[(ii)] $\HH^3(W,V)$ is a free abelian group, supporting a mixed
Hodge structure with weights $\geq 3$.
\end{enumerate}
\end{lemma}
\begin{proof} Let $W\setminus V=D$, and let $S\subset D$ be the union of
the singular locus of $D$, together with all irreducible components of $D$
of codimension $\ge 2$ in $X$. Then $D\setminus S=\coprod_jD_j$ where
$D_j\subset X\setminus S$ are irreducible, non-singular divisors.

We first observe that $\HH^i(W,W\setminus S)=0$ for $i\leq 3$, since
$S\subset W$ has (complex) codimension $\ge 2$. This implies that
$\HH^i(W,V)\to \HH^i(W\setminus S,V)$ are isomorphisms for $i\leq 3$, 
Since $(W\setminus S)\setminus V=\coprod_jD_j$, we have Thom-Gysin
isomorphisms $\HH^i(W\setminus S,V)\cong \oplus_j\HH^{i-2}(D_j)(-1)$ for
all $i\ge 0$ (where the Tate twist $(-1)$ increases the weights
by 2). In particular, we have $\HH^i(W\setminus S,V)=0$ for $i<2$, 
$\HH^0(D_j)=\bbZ$ (the trivial MHS), and  
$\HH^1(D_j)={\rm Hom}\,(H_1(D_j,\bbZ),\Z)$ is a torsion-free abelian
group,
which supports a MHS of weights $\geq 1$. 
\end{proof}
\begin{cor}\label{corr}
$\HH^{1}(X) \isom \HH^{1}(Y)$.
\end{cor}
\begin{proof} We have a factorization $\HH^1(X)\to\HH^1(X')\to
\HH^1(Y')=\HH^1(Y)$, since $Y=Y'\subset X'\subset X$.  By
Theorem~\ref{smt}, we have that $\HH^1(X')\to \HH^1(Y')=\HH^1(Y)$ is an
isomorphism. In particular, $\HH^1(X')$ supports a pure Hodge structure of
weight 1. Consider the exact cohomology sequence
\[\HH^1(X,X')\to\HH^1(X)\to\HH^1(X')\to\HH^2(X,X')\to \cdots\]
By Lemma~\ref{tf} with $W=X$, $V=X'$, we have $\HH^1(X,X')=0$, while 
$\HH^{2}(X,X')$ is torsion free, and it is pure of weight 2,
generated by the cohomology classes of the irreducible divisors in
$X\setminus X'$. Hence  the boundary map $ \HH^{1}(X') \to \HH^{2}(X,X')$
is the zero map. Thus $\HH^{1}(X) \to \HH^{1}(X')$ is an isomorphism.  
\end{proof}
\begin{cor}\label{Pic^0cor}
$\Pic^{0}(X) \to \Pic^{0}(Y)$ is an isomorphism.
\end{cor}
\begin{proof} The isomorphism in Corollary~\ref{corr} is compatible with
the respective Hodge structures, and so by Lemma~\ref{red-to-top}, induces
an isomorphism on $\Pic^0$ groups.
\end{proof}

\begin{prop}
Let $X$ and $Y$ be as in theorem \ref{GL2}. One then has an exact sequence
$$0 \to K \to \HH^{2}(X) \to \HH^{2}(Y) \to Q \to 0$$ 
where $K$ (as in Theorem \ref{GL2}) is generated by divisors which map to
points under $\pi$ and $Q$ is generated by the divisors in $Y$ which map
to points under $\pi$.
\end{prop}
\begin{proof}
Consider the diagram:
\begin{equation}\label{prime}
\begin{diagram}{cccccc}
\HH^1(X)&\by{a_1}\HH^1(X'')\by{a_2} &\HH^2(X,X'') &\by{a_3} \HH^2(X) \by{a_4} &
\HH^2(X'') &\by{a_5}\HH^3(X,X'')\\  
\quad\downarrow\psi_1&\quad\downarrow\psi_2&\quad\downarrow\psi_3&
\quad\downarrow\psi_4&\quad\downarrow \psi_5&\quad\downarrow\\
\HH^1(Y)&\by{b_1}\HH^1(Y'')\by{b_2} &\HH^2(Y,Y'') &\by{b_3}
\HH^2(Y) \by{b_4} &
\HH^2(Y'') &\by{b_5}\HH^3(Y,Y'')\\
\end{diagram}
\end{equation}
Here the horizontal sequences are the cohomology long exact
sequences corresponding to suitable pairs. The Proposition amounts to the
assertions  that there are isomorphisms
\[\ker \psi_3\by{\cong}\ker \psi_4,\]
\[\coker\psi_3\by{\cong}\coker\psi_4.\]

From Theorem~\ref{smt}, $\psi_2$ and $\psi_5$ are isomorphisms, while
$\psi_1$ is an isomorphism from Corollary~\ref{corr}. We claim that
$\psi_5$ induces an isomorphism 
\[{\rm image}\,a_4\by{\cong}{\rm image}\,b_4.\]
This follows by an argument using weights. Let
$W_2\HH^2(X'')\subset \HH^2(X'')$,
$W_2\HH^2(Y'')\subset \HH^2(Y'')$ be the subgroups
obtained as inverse images of the corresponding weight
subspaces of cohomology with rational coefficients. 
Since $\psi_5$ is an isomorphism of mixed Hodge
structures, it induces an isomorphism
$W_2\HH^2(X'')\cong W_2\HH^2(Y'')$. By Lemma~\ref{tf},
$\HH^3(X,X'')$ and $\HH^3(Y,Y'')$ are torsion free,
and have  weights $\geq 3$, while $\HH^2(X)$,
$\HH^2(Y)$ are pure of weight 2. Hence we have
\[W_2\HH^2(X'')={\rm image}\,a_4\]
\[W_2\HH^2(Y'')={\rm image}\,b_4,\]
and so $\psi_5$ induces an isomorphism between these
image subgroups.

Thus we have a commutative diagram with exact rows, and vertical
isomorphisms as shown.
\[\begin{diagram}{cccccc}
\HH^1(X)&\by{a_1}\HH^1(X'')\by{a_2} &\HH^2(X,X'') &\by{a_3} \HH^2(X) \by{a_4} &
{\rm image}\,a_4 &\to 0\\  
\cong\downarrow\psi_1&\cong\downarrow\psi_2&\quad\downarrow\psi_3&
\quad\downarrow\psi_4&\cong\downarrow\psi_5\\
\HH^1(Y)&\by{b_1}\HH^1(Y'')\by{b_2} &\HH^2(Y,Y'') &\by{b_3}
\HH^2(Y) \by{b_4} &
{\rm image}\,b_4 &\to 0\\
\end{diagram}
\]
A version of the 5-lemma  now implies that 
$\ker\psi_3\to\ker\psi_4$  and $\coker\psi_3\to\coker\psi_4$ are
isomorphisms, as desired. 
\end{proof}
\begin{cor}\label{NScor}
There is an exact sequence 
\[0\to K\to {\rm NS}\,(X)\to {\rm NS}\,(Y)\to Q\to 0.\]
\end{cor}
\begin{proof} This follows from the Proposition, and
Lemma~\ref{red-to-top}, since the explicit descriptions of $K$
and $Q$ imply that $K\subset {\rm NS}\,(X)$, and ${\rm NS}(Y)\onto Q$.
\end{proof}

\end{section}
\providecommand{\bysame}{\leavevmode\hbox to3em{\hrulefill}\thinspace}

\end{document}